\documentclass[12pt]{article}

\usepackage{amsmath,amssymb,amsthm,xspace,amscd}

\theoremstyle{remark}
\newtheorem{remark}{Remark}

\theoremstyle{plain}
\newtheorem{theorem}{Theorem}

\newtheorem{proposition}{Proposition}
\newtheorem{lemma}{Lemma}

\newcommand{\C}{{\mathbb C}}

\newcommand{\DR}{{\mathrm D}}
\newcommand{\Cr}{{\mathrm{Cr}}}
\newcommand{\Rat}{{\mathrm R}}
\newcommand{\Dist}{{\mathrm{Dist}}}
\newcommand{\Obig}{{\mathcal O}}
\newcommand{\Sch}{{\mathcal S}}
\newcommand{\diam}{{\mathrm{diam}}}

\title{On Cross-Ratio Distortion and Schwarz Derivative}

\author{
    A.~Teplinsky\thanks{Institute of Mathematics, Kiev, Ukraine}
    }

\begin{document}

\maketitle

\begin{abstract}
The asymptotical estimates for the cross-ratio distortion with respect to a smooth monotone function of one
variable in terms of its Schwarz derivative are established.
\end{abstract}

\section{Introduction, Definitions and Results}

Though the concept of cross-ratio of four consecutively connected segments has its origin in elementary geometry,
the question about good estimates on cross-ratio distortion with respect to a smooth function of one variable
arose in connection with studies in one-dimensional dynamics. Such tools were developed and applied first to a
great success in \cite{Yo84} to the case of critical circle maps and in \cite{MS} to the case of unimodal interval
maps. They also played an important role in \cite{Sw} (although dated before \cite{MS}, it refers to both
\cite{Yo84} and a preprint version of \cite{MS}). However natural it would seem to apply such tools to circle
diffeomorphisms, it was not done before the very recent works~\cite{KT-jams,T-umzh}, where some of the classical
results of Herman's theory~\cite{H,Y,SK,KO} were re-proven and even strengthened due, in part, namely to a
thorough investigation of the asymptotic expansions for cross-ratio distortion. The aim of this short paper is to prove optimal asymptotic estimates for cross-ratio distortion for both smooth and holomorphic cases without referring to one-dimensional dynamics, just in the elementary calculus framework. (Moreover, the only tool from the calculus we use is the Taylor's formula with the remainder term in asymptotic form.)

Let us start with the definitions. It is more convenient to talk about ratios and cross-ratios of points rather
than segments.

The {\em ratio} of three pairwise distinct points $x_1, x_2, x_3$ is
$$
\Rat(x_1,x_2,x_3)=\frac{x_1-x_2}{x_2-x_3},
$$
and the {\em ratio distortion} of those points with respect to the function $f$ is
$$
\DR(x_1,x_2,x_3;f)=\frac{\Rat(f(x_1),f(x_2),f(x_3))}{\Rat(x_1,x_2,x_3)}=
\frac{f(x_1)-f(x_2)}{x_1-x_2}:\frac{f(x_2)-f(x_3)}{x_2-x_3}.
$$

The {\em cross-ratio} of four pairwise distinct points $x_1, x_2, x_3, x_4$ is
$$
\Cr(x_1,x_2,x_3,x_4)=\frac{(x_1-x_2)(x_3-x_4)}{(x_2-x_3)(x_4-x_1)},
$$
whereas the {\em cross-ratio distortion} of those points with respect to $f$ is
\begin{eqnarray*}
\Dist(x_1,x_2,x_3,x_4;f)=\frac{\Cr(f(x_1),f(x_2),f(x_3),f(x_4))}{\Cr(x_1,x_2,x_3,x_4)}\\
=\frac{f(x_1)-f(x_2)}{x_1-x_2}:\frac{f(x_2)-f(x_3)}{x_2-x_3}\cdot
\frac{f(x_3)-f(x_4)}{x_3-x_4}:\frac{f(x_4)-f(x_1)}{x_4-x_1}.
\end{eqnarray*}

If the function $f$ is differentiable and its first derivative does not have zeros, then both ratio and
cross-ratio distortions are defined for not pairwise distinct points as well. Namely, these distortions can be
defined as the appropriate limits, or just by formally substituting $f'(a)$ for $\frac{f(a)-f(a)}{a-a}$ in the
definitions above. It is obvious that either $x_1=x_3$ or $x_2=x_4$ implies $\Dist(x_1,x_2,x_3,x_4;f)=1$.

As we find, the leading terms in the asymptotic expansion for cross-ratio distortion are directly related to the expression called `Schwarz derivative' that manifests itself in many considerations of one-dimensional real and complex dynamics. The {\em Schwarz derivative}, or {\em Schwarzian}, of a three times differentiable function $f$ at a point $x$ is given by
$$
\Sch f(x)=\frac{f'''(x)}{f'(x)}-\frac{3}{2}\left(\frac{f''(x)}{f'(x)}\right)^2
$$
as soon as that $f'(x)\ne0$. The connection between cross-ratio distortion and Schwarzian becomes evident if one
considers the two well-known facts about linear-fractional functions (a.k.a. `Moebius transformations'): on one
hand, $f$ is fractional-linear on $[A,B]$ if and only if $\Sch f\equiv0$ on $[A,B]$; on the other, $f$ is
fractional-linear on $[A,B]$ if and only if the cross-ratio distortion of any four points from $[A,B]$ with
respect to $f$ is equal to~1. Thus both Schwarzian and cross-ratio distortion in a sense measure how far is the
function $f$ from being fractional-linear. This is similar to the relation between the second derivative, ratio
distortion and non-linearity of a function. (A review of elementary facts known about cross-ratios and
Schwarzians can be found in \cite{MS}.)

Now we are ready to formulate our results. They are presented in a
series of four estimates related to different degrees of
smoothness: the first one applies to the case of smoothness $C^2$
and higher, the second one to $C^3$ and higher, the third one to
$C^4$ and higher, and the last one to the holomorphic case. Let us
remind that a domain $\Omega\subset\C$ is called {\em quasiconvex}
if there exists a constant $\Lambda\ge1$ such that for any two
points $a,b\in\Omega$ there exists a simple curve connecting them
such that its length does not exceed $\Lambda|a-b|$.

Note, that all the implicit constants, which are presented
throughout this paper in the form of $\Obig(\cdot)$, depend on the
function $f$ and its segment of definition $[A,B]$ only (in the
smooth case) or on the function $F$ and a chosen compact subset of
its domain of definition $\Omega$ only (in the holomorphic case).
For a (finite) set $M$, by $\diam M$ we denote its diameter, i.e.
the greatest distance between its points.

\begin{theorem}
Let $f\in C^r([A,B])$, and $f'$ does not have zeroes on $[A,B]$. Consider four arbitrary points $x_1, x_2, x_3,
x_4\in[A,B]$ and denote $\Delta=\diam\{x_1,x_2,x_3,x_4\}$. The stated below asymptotic estimates hold true.

In the case of $r=2+\alpha$, $\alpha\in[0,1]$, $\Delta\ne0$:
\begin{equation}\label{eq:th1}
\Dist(x_1,x_2,x_3,x_4;f)=1+(x_1-x_3)(x_2-x_4)\Obig(\Delta^{\alpha-1}).
\end{equation}

In the case of $r=3+\beta$, $\beta\in[0,1]$:
\begin{equation}\label{eq:th2}
\Dist(x_1,x_2,x_3,x_4;f)=1+(x_1-x_3)(x_2-x_4)\left(\frac{1}{6}\Sch
f(\theta)+\Obig(\Delta^{\beta})\right)
\end{equation}
with arbitrary
$\theta\in[\min\{x_1,x_2,x_3,x_4\},\max\{x_1,x_2,x_3,x_4\}]$.

In the case of $r=4+\gamma$, $\gamma\in[0,1]$:
\begin{eqnarray}
\Dist(x_1,x_2,x_3,x_4;f)=1+(x_1-x_3)(x_2-x_4)\nonumber\\
\times\left(\frac{1}{24}\bigl(\Sch f(x_1)+\Sch f(x_2)+\Sch f(x_3)+\Sch
f(x_4)\bigr)+\Obig(\Delta^{1+\gamma})\right).\label{eq:th3}
\end{eqnarray}

Let $F$ be a holomorphic function defined on a quasiconvex domain $\Omega\subset\C$ such that $F'$ does not have
zeroes in $\Omega$. Uniformly on compact subsets of $\Omega$, the following asymptotic estimate holds true:
\begin{eqnarray}
\Dist(z_1,z_2,z_3,z_4;F)=1+(z_1-z_3)(z_2-z_4)\nonumber\\
\times\left(\frac{1}{24}\bigl(\Sch F(z_1)+\Sch F(z_2)+\Sch F(z_3)+\Sch
F(z_4)\bigr)+\Obig(\Delta^2)\right),\label{eq:th4}
\end{eqnarray}
where in this case $\Delta=\diam\{z_1,z_2,z_3,z_4\}$.
\end{theorem}

\begin{remark}
We wish to stress it straight away that the leading terms in this
asymptotic expansion are not too hard to derive by themselves, whereas the
proof that the remainder term for $f\in C^r$ is
$(x_1-x_3)(x_2-x_4)\Obig(\Delta^{r-3})$ rather than just
$\Obig(\Delta^{r-1})$ is far from obvious (and it is clear that
the distances $|x_1-x_3|$ and $|x_2-x_4|$ can be much smaller than
$\Delta$). A similar remark applies to the holomorphic case.
\end{remark}

\section{Proof of Theorem~1}

Here we will consider the case $f\in C^{4+\gamma}([A,B])$, $\gamma\in[0,1]$, and prove the estimate
(\ref{eq:th3}). As it will become evident, the proofs of (\ref{eq:th1}), (\ref{eq:th2}) and (\ref{eq:th4}) follow
the same lines with very slight modifications.

Let us introduce notations $\phi_k=\frac{f^{(k+1)}(\theta)}{(k+1)!f'(\theta)}$ and $d_i=x_i-\theta$. Let
$x_1,x_2,\theta$ be arbitrary points from the segment $[A,B]$. It is easy to derive from the Taylor's expansions
for $f(x_1)$ ³ $f(x_2)$ with respect to the reference point $\theta$ that
\begin{equation}\label{eq:frac}
\frac{f(x_1)-f(x_2)}{f'(\theta)(x_1-x_2)}=1+P_1+P_2+P_3+\Obig\bigl((\diam\{x_1,x_2,\theta\})^{3+\gamma}\bigr),
\end{equation}
where $P_k=\phi_k\frac{d_1^{k+1}-d_2^{k+1}}{x_1-x_2}=\phi_k\sum_{j=0}^kd_1^jd_2^{k-j}$, $k\in\{1,2,3\}$, are the
symmetric polynomials of degree $k$ with respect to $d_1$ and $d_2$.

Before we start the actual proof, let us show a way that produces the leading terms of the asymptotic expansion straight away, although does not give the optimal estimate. Using the expansion
$\ln(1+t)=t-\frac{t^2}{2}+\frac{t^3}{3}+\Obig(t^4)$, we achieve
\begin{eqnarray}
\ln\frac{f(x_1)-f(x_2)}{f'(\theta)(x_1-x_2)}=P_1+\left[P_2-\frac{1}{2}P_1^2\right]\nonumber\\
+\left[P_3-P_1P_2+\frac{1}{3}P_1^3\right]
+\Obig\bigl((\diam\{x_1,x_2,\theta\})^{3+\gamma}\bigr)\label{eq:lnfrac}
\end{eqnarray}
(here and in what follows, in square brackets we group up terms of the same order). Now, if one would simply
calculate $\ln\Dist(x_1,x_2,x_3,x_4;f)$ as the sum of the four expressions
$$
\ln\frac{f(x_1)-f(x_2)}{f'(\theta)(x_1-x_2)}
-\ln\frac{f(x_2)-f(x_3)}{f'(\theta)(x_2-x_3)}+\ln\frac{f(x_3)-f(x_4)}{f'(\theta)(x_3-x_4)}
-\ln\frac{f(x_4)-f(x_1)}{f'(\theta)(x_4-x_1)},
$$
substituting the corresponding variants of the expansion (\ref{eq:lnfrac}), then after appropriate
transformations the formula
\begin{eqnarray*}
\ln\Dist(x_1,x_2,x_3,x_4;f)=(x_1-x_3)(x_2-x_4)\\
\times\frac{1}{24}\bigl(\Sch f(x_1)+\Sch f(x_2)+\Sch f(x_3)+\Sch f(x_4)\bigr)+\Obig(\Delta^{3+\gamma})
\end{eqnarray*}
will be obtained. However, the remainder term in it is not what we are looking for. The optimal estimate
(\ref{eq:th3}) cannot be proven in such a direct way, and so we shall take a roundabout path in order to extract the multiple $(x_1-x_3)(x_2-x_4)$ from that remainder term.

\begin{lemma}\label{lemma:equiv}
The following exact equalities take place:
\begin{eqnarray}
(x_2-x_3)\bigl(\DR(x_1,x_2,x_3;f)-1\bigr)=(x_1-x_3)\bigl(\DR(x_2,x_1,x_3;f)-1\bigr)\DR(x_1,x_3,x_2;f);\label{eq:equiv1}\\
(x_2-x_3)(x_1-x_4)\bigl(\Dist(x_1,x_2,x_3,x_4;f)-1\bigr)\nonumber\\
=(x_1-x_3)(x_2-x_4)\bigl(\Dist(x_2,x_1,x_3,x_4;f)-1\bigr)\Dist(x_1,x_3,x_2,x_4;f).\label{eq:equiv2}
\end{eqnarray}
\end{lemma}
{\em Proof.}\/ We will prove both (\ref{eq:equiv1}) and (\ref{eq:equiv2}) under the condition that
$x_1,x_2,x_3,x_4$ are pairwise distinct. The cases, when some of those points coincide, are very easy to check
directly or can be reached from the pairwise distinct case by appropriate limit transitions.

One can see that $\Rat(x_3,x_1,x_2)+\Rat(x_3,x_2,x_1)=-1$;
\\
also $\Rat(f(x_3),f(x_1),f(x_2))+\Rat(f(x_3),f(x_2),f(x_1))=-1$. Hence,
$$
\frac{x_2-x_3}{x_1-x_2}-\frac{f(x_2)-f(x_3)}{f(x_1)-f(x_2)}=-\frac{x_1-x_3}{x_2-x_1}+\frac{f(x_1)-f(x_3)}{f(x_2)-f(x_1)},
$$
which implies
$$
(f(x_2)-f(x_3))(\DR(x_1,x_2,x_3;f)-1)=(f(x_1)-f(x_3))(\DR(x_2,x_1,x_3;f)-1).
$$
The latter formula is easily transformed into the equality (\ref{eq:equiv1}).

Since $(x_2-x_3)(x_4-x_1)-(x_1-x_3)(x_4-x_2)=(x_1-x_2)(x_3-x_4)$, we have
$\Cr(x_2,x_3,x_4,x_1)+\Cr(x_1,x_3,x_4,x_2)=1$; also
$\Cr(f(x_2),f(x_3),f(x_4),f(x_1))+\Cr(f(x_1),f(x_3),f(x_4),f(x_2))=1$. Hence,
\begin{eqnarray*}
\frac{(x_2-x_3)(x_4-x_1)}{(x_1-x_2)(x_3-x_4)}-\frac{(f(x_2)-f(x_3))(f(x_4)-f(x_1))}{(f(x_1)-f(x_2))(f(x_3)-f(x_4))}\\
=-\frac{(x_1-x_3)(x_4-x_2)}{(x_2-x_1)(x_3-x_4)}+\frac{(f(x_1)-f(x_3))(f(x_4)-f(x_2))}{(f(x_2)-f(x_1))(f(x_3)-f(x_4))},
\end{eqnarray*}
and therefore
\begin{eqnarray*}
(f(x_2)-f(x_3))(f(x_4)-f(x_1))(\Dist(x_1,x_2,x_3,x_4;f)-1)\\
=(f(x_1)-f(x_3))(f(x_4)-f(x_2))(\Dist(x_2,x_1,x_3,x_4;f)-1),
\end{eqnarray*}
which is easy to transform into the equality (\ref{eq:equiv2}).

Lemma~\ref{lemma:equiv} is proven.

Consider the expression
\begin{eqnarray*}
 Q(\theta,x_1,x_2,x_3)=\phi_1+\bigl[\phi_2(d_1+d_2+d_3)-\phi_1^2(d_2+d_3)\bigr]\\
+\bigl[\phi_3(d_1^2+d_2^2+d_3^2+d_1d_2+d_2d_3+d_3d_1)\\
-\phi_1\phi_2((d_2^2+d_2d_3+d_3^2)+(d_2+d_3)(d_1+d_2+d_3))+\phi_1^3(d_2+d_3)^2\bigr],
\end{eqnarray*}
which in the sequel we will denote simply as $Q_{123}$.

\begin{proposition}\label{prop:th3forD}
Let $f\in C^{4+\gamma}([A,B])$, $\gamma\in[0,1]$, and $f'>0$. For any four points $x_1, x_2, x_3,\theta\in[A,B]$
the following asymptotic estimate takes place:
\begin{equation}\label{eq:th3forD}
\DR(x_1,x_2,x_3;f)=1+(x_1-x_3)\bigl(Q_{123}+\Obig(\Delta_\theta^{2+\gamma})\bigr),
\end{equation}
where $\Delta_\theta=\diam\{x_1,x_2,x_3,\theta\}$.
\end{proposition}

\begin{remark}
An arbitrary choice of $\theta$ in Proposition~\ref{prop:th3forD} makes that form of the asymptotic estimate the most
general, giving an opportunity to produce different variants of the estimate (\ref{eq:th3forD}) for different
specific $\theta$ (in particular, one can consider the variants with $\theta=x_1$, $\theta=x_2$ or $\theta=x_3$).
\end{remark}

First, let us prove the following lemma concerning the dependence of $Q_{123}$ on $\theta$.
\begin{lemma}\label{lemma:tildeQ-Q}
Let $x_1,x_2,x_3,\theta,\tilde\theta\in[A,B]$, and $\tilde Q_{123}=Q(\tilde\theta,x_1,x_2,x_3)$. The following
asymptotic estimate takes place: $\tilde Q_{123}-Q_{123}=\Obig(|\delta|^{2+\gamma})$, where
$\delta=\tilde\theta-\theta$.
\end{lemma}
{\em Proof.}\/ Let us find the partial asymptotic expansions for
$\tilde\phi_k=\frac{f^{(k+1)}(\tilde\theta)}{(k+1)!f'(\tilde\theta)}$ in terms of $\phi_k$ with respect to the
powers of $\delta$. In the case of $k=1$ we write
\begin{equation}\label{tildephi1}
\tilde\phi_1=\frac{1}{2}\frac{f''(\tilde\theta)/f'(\theta)}{f'(\tilde\theta)/f'(\theta)}
=\frac{\phi_1+3\phi_2\delta+6\phi_3\delta^2+\Obig(|\delta|^{2+\gamma})}
{1+2\phi_1\delta+3\phi_2\delta^2+\Obig(|\delta|^3)},
\end{equation}
which implies (in view of the expansion
$\frac{1}{1+t}=1-t+t^2+\Obig(t^3)$ and after noticing that the
absolute value of the denominator in (\ref{tildephi1}) is confined
between two positive constants)
$$
\tilde\phi_1=\phi_1+[3\phi_2-2\phi_1^2]\delta+[6\phi_3-9\phi_2\phi_1+4\phi_1^3]\delta^2+\Obig(|\delta|^{2+\gamma}).
$$
Similarly obtain
$$
\tilde\phi_2=\frac{1}{2}\frac{f'''(\tilde\theta)/f'(\theta)}{f'(\tilde\theta)/f'(\theta)}
=\frac{\phi_2+4\phi_3\delta+\Obig(|\delta|^{1+\gamma})}{1+2\phi_1\delta+\Obig(|\delta|^2)}
=\phi_2+[4\phi_3-2\phi_3\phi_2]\delta+\Obig(|\delta|^{1+\gamma})
$$
and, finally, $\tilde\phi_3=\phi_3+\Obig(|\delta|^{\gamma})$.

Now, substitute the derived expressions together with $\tilde d_i=x_i-\tilde\theta=d_i-\delta$, $i\in\{1,2,3\}$,
into $\tilde Q_{123}$, subtract $Q_{123}$, and after transformations get the estimate of the lemma.
Lemma~\ref{lemma:tildeQ-Q} is proven.

{\em Proof of Proposition~\ref{prop:th3forD}.}\/ According to Lemma~\ref{lemma:tildeQ-Q}, it is enough to prove
the estimate (\ref{eq:th3forD}) for any single point $\theta\in[\min\{x_1, x_2, x_3\},\max\{x_1, x_2, x_3\}]$, and
that will imply that (\ref{eq:th3forD}) is true for each $\theta\in[A,B]$. However, we will not specify the
choice of $\theta$ in this proof, imposing only the condition $\theta\in[\min\{x_1, x_2, x_3\},\max\{x_1, x_2,
x_3\}]$. (A constructivist reader is welcome to assume $\theta=x_1$, although that will not simplify the
expressions.) This condition implies $\Delta_\theta=\diam\{x_1,x_2,x_3\}$, which we will denote by $\Delta_{123}$
during this proof.

It follows from the definition of ratio distortion that
\begin{equation}\label{eq:A-B}
\DR(x_1,x_2,x_3;f)=1+\frac{c_{12}-c_{23}}{1+c_{23}},
\end{equation}
where $c_{12}=\frac{f(x_1)-f(x_2)}{f'(\theta)(x_1-x_2)}-1$, $c_{23}=\frac{f(x_2)-f(x_3)}{f'(\theta)(x_2-x_3)}-1$.
According to (\ref{eq:frac}), we have
$$
c_{12}=\phi_1(d_1+d_2)+\phi_2(d_1^2+d_1d_2+d_2^2)+\phi_3(d_1^3+d_1^2d_2+d_1d_2^2+d_2^3)+\Obig(\Delta_{123}^{3+\gamma}),
$$
$$
c_{23}=\phi_1(d_2+d_3)+\phi_2(d_2^2+d_2d_3+d_3^2)+\phi_3(d_2^3+d_2^2d_3+d_2d_3^2+d_3^3)+\Obig(\Delta_{123}^{3+\gamma}).
$$
Substitute these expressions into (\ref{eq:A-B}) in view of
$\frac{1}{1+t}=1-t+t^2+\Obig(t^3)$ (noticing that the absolute
value of the denominator $1+c_{23}$ is confined between two
positive constants again) and after transformations get
\begin{equation}\label{eq:th3forDbad}
\DR(x_1,x_2,x_3;f)=1+(x_1-x_3)Q_{123}+\Obig(\Delta_{123}^{3+\gamma}).
\end{equation}
%$\frac{\min|f'|}{\max|f'|}$ and $\frac{\max|f'|}{\min|f'|}$

The estimate (\ref{eq:th3forDbad}) implies (\ref{eq:th3forD}) in the case when the points $\theta$ and $x_2$ lie
between the points $x_1$ and $x_3$ (so that $\Delta_\theta=\Delta_{123}=|x_1-x_3|$). Thus, in that case the lemma
is proven.

Now suppose that $\theta$ and $x_1$ lie between $x_2$ and $x_3$, so that $\Delta_\theta=\Delta_{123}=|x_2-x_3|$.
Having transposed the points in (\ref{eq:th3forDbad}) as necessary, we obtain
$$
\DR(x_2,x_1,x_3;f)=1+(x_2-x_3)\bigl(Q_{213}+\Obig(\Delta_{123}^{2+\gamma})\bigr),
$$
$$
\DR(x_1,x_3,x_2;f)=1+(x_1-x_2)Q_{132}+\Obig(\Delta_{123}^{3+\gamma}),
$$
where $Q_{213}$ and $Q_{132}$ are obtained of $Q_{123}$ by corresponding transpositions of variables $d_1$, $d_2$
and $d_3$. Using the equality (\ref{eq:equiv1}), we get
\begin{eqnarray}
\DR(x_1,x_2,x_3;f)=1+(x_1-x_3)\nonumber\\
\times\bigl(Q_{213}+\Obig(\Delta_{123}^{2+\gamma})\bigr)
\bigl(1+(d_1-d_2)Q_{132}+\Obig(\Delta_{123}^{3+\gamma})\bigr).\label{eq:th3forDinter}
\end{eqnarray}
It is easy to calculate that
$$
Q_{123}-Q_{213}=(d_1-d_2)(\phi_1^2+[2\phi_2\phi_1(d_1+d_2+d_3)-\phi_1^3(d_1+d_2+2d_3)]),
$$
$$
Q_{213}Q_{132}=\phi_1^2+[2\phi_2\phi_1(d_1+d_2+d_3)-\phi_1^3(d_1+d_2+2d_3)]+\Obig(\Delta_{123}^2),
$$
so (\ref{eq:th3forDinter}) implies (\ref{eq:th3forD}) indeed.

The case when $\theta$ and $x_3$ lie between $x_1$ and $x_2$, is done similarly. Proposition~\ref{prop:th3forD}
is proven.

{\em Proof of (\ref{eq:th3}).}\/ Let $\theta\in[\min\{x_1, x_2, x_3, x_4\},\max\{x_1, x_2, x_3, x_4\}]$. Using the
definitions of $\DR$ and $\Dist$ and Proposition~\ref{prop:th3forD}, we get
\begin{eqnarray*}
\Dist(x_1,x_2,x_3,x_4;f)=\DR(x_1,x_2,x_3;f)\cdot\DR(x_3,x_4,x_1;f)\\
=\bigl(1+(x_1-x_3)(S_{123}+\Obig(\Delta^{2+\gamma}))\bigr)
\bigl(1+(x_3-x_1)(S_{341}+\Obig(\Delta^{2+\gamma}))\bigr)\\
=1+(x_1-x_3)\bigl(S_{123}-S_{341}-(x_1-x_3)S_{123}S_{341}+\Obig(\Delta^{2+\gamma})\bigr).
\end{eqnarray*}
Simple transformations show that
\begin{eqnarray*}
 S_{123}-S_{341}-(d_1-d_3)S_{123}S_{341}\\
=(d_2-d_4)\bigl((\phi_2-\phi_1^2)+(\phi_3-2\phi_2\phi_1+\phi_1^3)(d_1+d_2+d_3+d_4)\bigr)+\Obig(\Delta^3).
\end{eqnarray*}
It is time to notice that $\phi_2-\phi_1^2=\frac{1}{6}\Sch
f(\theta)$, $\phi_3-2\phi_2\phi_1+\phi_1^3=\frac{1}{24}(\Sch
f)'(\theta)$, and $\Sch f(\theta)+(\Sch f)'(\theta)d_i=\Sch
f(x_i)+\Obig(|d_i|^{1+\gamma})$ for $i\in\{1,2,3,4\}$, so that we
finally obtain
\begin{eqnarray}
\Dist(x_1,x_2,x_3,x_4;f)=1+(x_1-x_3)\nonumber\\
\times\left((x_2-x_4)\frac{1}{24}\sum_{i=1}^4\Sch f(x_i)+\Obig(\Delta^{2+\gamma})\right).\label{eq:th3bad}
\end{eqnarray}
The role of (\ref{eq:th3bad}) in this proof is similar to the role of (\ref{eq:th3forDbad}) in the proof of
Proposition~\ref{prop:th3forD}. Namely, in the case when $x_1$ and $x_3$ lie between $x_2$ and $x_4$ we have
$\Delta=|x_2-x_4|$, and hence (\ref{eq:th3bad}) implies (\ref{eq:th3}). Thus, in that case the theorem is proven.
Notice, that if $x_2$ and $x_4$ lie between $x_1$ and $x_3$, then the theorem is proven as well due to the
symmetry $\Dist(x_1,x_2,x_3,x_4;f)=\Dist(x_2,x_1,x_4,x_3;f)$.
%(it is enough to swap $x_1$ with $x_2$ and $x_3$with $x_4$ in this proof).

Now suppose that $x_2$ and $x_3$ lie between $x_1$ and $x_4$, so that $\Delta=|x_1-x_4|$. Obvious transpositions
of points in (\ref{eq:th3bad}) lead to
$$
\Dist(x_2,x_1,x_3,x_4;f)=1+(x_2-x_3)(x_1-x_4)\left(\frac{1}{24}\sum_{i=1}^4\Sch
f(x_i)+\Obig(\Delta^{1+\gamma})\right),
$$
$$
\Dist(x_1,x_3,x_2,x_4;f)=1+\Obig(\Delta^{2}),
$$
and (\ref{eq:th3}) follows from the equality (\ref{eq:equiv2}). Thus the theorem is proven in this case, too. By
symmetry, it is proven also for the case when $x_1$ and $x_4$ lie between $x_2$ and $x_3$.

Finally, the case of $x_1$ and $x_2$ lying between $x_3$ and $x_4$ (and the symmetric one, with $x_3$ and $x_4$
between $x_1$ and $x_2$) is considered similarly. Thus (\ref{eq:th3}) is proven.

It is quite obvious now that the proofs of (\ref{eq:th1}) and (\ref{eq:th2}) are easily obtained from the proof
of (\ref{eq:th3}) by cutting off all the derived partial asymptotic expansions at appropriate lower-order terms.

It is also not hard to check that (\ref{eq:th4}) is proven by
following the lines of the proof of (\ref{eq:th3}) with $\gamma=1$
in appropriate settings. All the statements of the form ``$a$ lies
between $b$ and $c$'' are to be replaced with
``$\diam\{a,b,c\}=|b-c|$'', whereas for ``$b\in[\min M ,\max M]$''
for a finite set $M$ one has to substitute ``$\diam\bigl(\{b\}\cup
M\bigr)=\diam M$''.

Theorem~1 is proven.

{\em Acknowledgements.}
The author is grateful to Konstantin Khanin, Ilia Binder and Welington de Melo for
valuable comments.

\end{document}